\documentclass[16pt]{article}

\usepackage{graphicx}
\usepackage{graphpap}
\usepackage[mathscr]{eucal}
\usepackage{amsmath}
\usepackage{amsfonts}
\usepackage{amssymb}
\usepackage{amsthm}
\usepackage{amstext}
\usepackage{amscd}
\usepackage{makeidx}
\usepackage{graphics}

 \newtheorem{dfn}{Definition}
 \newtheorem{thm}{Theorem}
 \newtheorem{lem}{Lemma}
 \newtheorem{cor}{Corollary}

 \newtheorem{op}{Open problem}
 \newtheorem{con}{Conjecture}
 \newcommand{\bi}{\begin{itemize}}
 \newcommand{\ei}{\end{itemize}}
 \newcommand{\bd}{\begin{dfn}}
 \newcommand{\ed}{\end{dfn}}
 \newcommand{\bt}{\begin{thm}}
 \newcommand{\et}{\end{thm}}
 \newcommand{\bl}{\begin{lem}}
 \newcommand{\el}{\end{lem}}
 \newcommand{\bc}{\begin{cor}}
 \newcommand{\ec}{\end{cor}}
 \newcommand{\bcon}{\begin{con}}
 \newcommand{\econ}{\end{con}}
 \newcommand{\bop}{\begin{op}}
 \newcommand{\eop}{\end{op}}
 \newcommand{\bea}{\begin{eqnarray}}
 \newcommand{\eea}{\end{eqnarray}}
 \newcommand{\ebox}{\hfill $\Box$\\\vspace{0.2cm}}
 \newcommand{\pr}{\noindent{\bf Proof.}\ }
 \newcommand{\noi}{\noindent}

\begin{document}

\title{Magic labelings of distance at most 2}
\author{
\normalsize Rinovia Simanjuntak\thanks{This article is written as a class project of MA6151 Topics in Discrete Mathematics: Introduction to Graph Labeling, 2011/2012, Mathematics Master Program, Institut Teknologi Bandung}, Mona Elviyenti, Mohammad Nafie Jauhari,\\
\normalsize Alfan Sukmana Praja, and Ira Apni Purwasih\\
\vspace{-1mm}
\normalsize Combinatorial Mathematics Research Group\\
\vspace{-1mm}
\normalsize Faculty of Mathematics and Natural Sciences\\
\vspace{-1mm}
\normalsize Institut Teknologi Bandung, Bandung 40132, Indonesia\\
\vspace{-2mm}
\scriptsize {\small {\bf e-mail}: {\tt rino@math.itb.ac.id}}\\\\
}
\date{}
\maketitle

\vspace*{-4mm}
\begin{abstract}
For an arbitrary set of distances $D\subseteq \{0,1, \ldots, d\}$, a graph $G$ is said to be $D$-distance magic if there exists a bijection $f:V\rightarrow \{1,2, \ldots , v\}$ and a constant {\sf k} such that for any vertex $x$, $\sum_{y\in N_D(x)} f(y) ={\sf k}$, where $N_D(x) = \{y \in V| d(x,y) \in D\}$.

In this paper we study some necessary or sufficient conditions for the existence of $D$-distance magic graphs, some of which are generalization of conditions for the existence of $\{1\}$-distance magic graphs. More specifically, we study $D$-distance magic labelings for cycles and $D$-distance magic graphs for $D\subseteq\{0,1,2\}$.
\end{abstract}

\section{Introduction}

\noi As standard notation, assume that $G$=$G(V,E)$ is a finite, simple, and undirected graph with $v$ vertices, $e$ edges, and diameter $d$. By a {\em labeling} we mean a one-to-one mapping that carries a set of graph elements onto a set of
numbers, called {\em labels}.

\noi The notion of distance magic labeling was introduced separately in the PhD thesis of Vilfred \cite{Vi94}
in 1994 and an article by Miller {\it et. al} \cite{MRS03} in 2003. A {\em distance magic labeling}, or $\Sigma$ labeling, is a bijection $f:V\rightarrow \{1,2, \ldots , v\}$ with the property that there is a constant {\sf k} such that at any vertex $x$,
$\sum_{y\in N(x)} f(y) ={\sf k}$, where $N(x)$ is the open neighborhood of $x$, i.e., the set of vertices adjacent to $x$. This labeling was introduced due to two different motivations; as a tool in utilizing magic squares into graphs and as a natural extension of previously known graph labelings: magic labeling \cite{Sed64,KR70} and radio labeling (which is distance-based) \cite{GY92}.

\noi In the last decade, many results on distance magic labeling have been published. Several families of graphs have been showed to admit the labeling, for instance circulant graphs \cite{CF}, bipartite graphs \cite{MRS03,ARSP04,Be09,CG}, tripartite graphs \cite{MRS03}, regular multipartite graphs \cite{Vi94,MRS03}, Cartesian product graphs \cite{Ji99,Ra08}, lexicographic product graphs \cite{MRS03,SAS09,AC,ACPT}, and joint product graphs \cite{SAA09}. Constructions of distance magic graphs have also been studied: construction producing regular graphs was studied in \cite{Fr07,FKK06,FKK11,KFK12} and non-regular graphs in \cite{SFMRW09,KS12}; the constructions utilize Kotzig array and magic rectangle. 

\noi It was proved in \cite{Vi94} that every graph is a subgraph of a distance magic graph. A stronger result that every graph is an induced subgraph of a regular distance magic graph was then proved in \cite{ARSP04}. A yet stronger result can also be found in \cite{RSP04} where it is stated that every graph $H$ is an induced subgraph of a Eulerian distance magic graph $G$ where the chromatic number of $H$ is the same as $G$. All these results showed that there is no forbidden subgraph characterization for distance magic graph. Additionally, an application of the labeling in designing incomplete tournament is introduced in \cite{FKK06}. For more results, please refer to a recent survey article in \cite{AFK13}.

\noi Jinah \cite{Ji99} introduced a variation of distance magic labeling. A $\Sigma ^{'}$ labeling, is a bijection $f:V\rightarrow \{1,2, \ldots , v\}$ with the property that there is a constant {\sf k} such that at any vertex $x$,
$\sum_{y\in N[x]} f(y) ={\sf k}$, where $N[x]$ is the closed neighborhood of $x$, i.e., the set containing $x$ and all vertices adjacent to $x$. It was stated in \cite{OS11b} that there does not exist a graph of even order that admits both distance magic and $\Sigma ^{'}$ labelings. As for graphs of odd order, the path $P_3$ is one example of a graph admitting both labelings. In the same article, it was also showed that a graph is distance magic if and only if its complement is $\Sigma ^{'}$-labeled.

\noi Recently O'Neal and Slater \cite{OS11b} generalized the notion of distance magic labeling to an arbitrary set of distances $D\subseteq \{0,1, \ldots, d\}$. As in distance magic labeling, the domain of the new labeling is the set of all vertices and the codomain is $\{1,2, \ldots , v\}$. We define the {\em $D$-weight} of each vertex $x$ in $G$, denoted by $w(x)$, to be the sum of labels of the vertices at distance $k$ to $x$, where $k \in D$. If all vertices in $G$ have the same weight, we call the labeling a {\em $D$-distance magic labeling}. More formally, we have the following definition.
\bd {\em \cite{OS11b}}
A bijection $f:V\rightarrow \{1,2, \ldots , v\}$ is said to be a {\bf $D$-distance magic labeling} if there exists a {\bf $D$-distance magic constant} {\sf k} such that for any vertex $x$, $w(x)=\sum_{y\in N_D(x)} f(y)={\sf k}$, where $N_D(x) = \{y \in V| d(x,y) \in D\}$. A graph admitting a $D$-distance magic labeling is called {\bf $D$-distance magic}.
\ed

\noi Clearly, a distance-magic labeling is a $\{1\}$-distance magic labeling and a $\Sigma ^{'}$ labeling is a $\{0,1\}$-distance magic labeling. Rewriting the results in \cite{OS11b}, we have the following relations between $\{1\}$-distance magic and $\{0,1\}$-distance magic labeling.
\bl \emph{\cite{OS11b}}
There does not exist a graph of even order that admits both $\{1\}$-distance magic and $\{0,1\}$-distance magic labelings.
\el
\bl \emph{\cite{OS11b}}
A graph is $\{1\}$-distance magic if and only if its complement is $\{0,1\}$-distance magic.
\el

\noi In this paper we study properties of $D$-distance magic labelings for a distance set $D$, where $D\subseteq \{0,1, \ldots, d\}$. Obviously, the only $\{0\}$-distance magic graph is the trivial graph, and so we exclude $D=\{0\}$ from our consideration. Additionally, we also study $D$-distance magic labelings for $D\subseteq\{0,1,2\}$.

\section{Some general results}

\noi In this section, we study some necessary and sufficient conditions for the existence of $D$-distance magic graphs for particular distance sets $D$, $D\subseteq \{0,1, \ldots, d\}$ and $D\neq \{0\}$. Unless stated, we shall exclude the trivial graph from consideration. We start by generalizing some properties of $\{1\}$-distance magic graphs presented in \cite{MRS03}.

\noi In \cite{MRS03} it was proved that there does not exist a $\{1\}$-distance magic labeling for $r$-regular graph with odd $r$. The next result generalize this idea to arbitrary neighborhood sets. Graph $G$ is defined to be {\em $(D,r)$-regular} if for all $v\in V (G), |N_D(v)|=r$, that is, all $D$-neighborhoods have the same cardinality.
\bl {\em \cite{OS11b}}
Let $G$ be a graph of even order. If $G$ is $(D,r)$-regular with odd $r$, then $G$ is not $D$-distance magic.
\label{(D,r)-regular}
\el

\noi Another result can be found in \cite{MRS03} is that if a graph $G$ contains two vertices $x$ and $y$ such that $|N(x) \cap N(y)|=
d(x)-1 = d(y)-1$ then $G$ is not $\{1\}$-distance magic. We shall generalize the idea to $D$-distance magic graphs.
\bl
If a graph $G$ contains two distinct vertices $x$ and $y$ such that $|N_D(x) \cap N_D(y)|= |N_D(x)|-1 = |N_D(y)|-1$ then $G$ is not $D$-distance magic.
\label{d-1}
\el
\begin{proof} Suppose $G$ is $D$-distance magic and let $x'$ ($y'$, respectively) be the one vertex in $N_D(x)-N_D(y)$ ($N_D(y)-N_D(x)$, respectively).
Then $\sum_{u\in N_D(x)} f(u) = w(x) = w(y) = \sum_{u\in N_D(y)} f(u)$, and so $f(x')=f(y')$, a contradiction.
\end{proof}

\noi The following two lemmas also give necessary conditions connected to the $D$-neighborhood of vertices in the graph.
\bl
If $G$ contains a vertex $x$ with $N_D(x)=\emptyset$ then $G$ is not $D$-distance magic.
\el
\begin{proof} Suppose $G$ is $D$-distance magic. Since $D \subseteq \{0,1, \ldots, d\}$ then there is a vertex $y$ where $N_D(y)\neq\emptyset$ and so $w(y)\neq 0$. However $w(x)= 0$, a contradiction.
\end{proof}

\bl
If $G$ contains two distinct vertices $x$ and $y$ such that $N_D(x)\subseteq N_D(y)$ then $G$ is not $D$-distance magic.
\el
\pr Suppose $G$ is $D$-distance magic. Since $w(x)=w(y)$, then $\sum_{u\in N_D(y)-N_D(x)} f(u)=0$, a contradiction. \ebox

\noi Properties of $D$-distance magic graphs can also be found in \cite{OS13}, the most important is the uniqueness of the $D$-distance magic constant.
\bd
A function $g: V(G)\rightarrow R^+ = [0,\infty)$ is said to be a {\bf $D$-neighborhood fractional dominating function} if for every $v\in V(G)$, $\sum_{u\in N_D(v)} g(u) \geq 1$. The {\bf $D$-neighborhood fractional domination number} of $G$, denoted by $\gamma_f(G;D)$, is defined as $\gamma_f(G;D) = \min \{\sum_{v\in V(G)} g(v)|g \textrm{ is a } $D$-\textrm{neighborhood fractional dominating function}\}$.
\ed
\bt {\em \cite{OS13}}
If graph $G$ is $D$-distance magic, then its $D$-distance magic constant ${\sf k} = \frac{n(n+1)}{2\gamma_f(G;D)}$.
\et

%
%
%
\noi The following two lemmas deal with existence of $D$-distance magic graphs for particular $D$.
\bl
If each vertex in $G$ has a unique vertex at distance $d$ then $G$ is $\{1,2,\ldots,d-1\}$-distance magic.
\label{1k-1}
\el
\begin{proof} 
We define a labeling $f$ such that if a vertex $x$ is labeled with $i$ then the unique vertex at distance $d$ from $x$ is labeled with $v+1-i$. Thus, for every vertex $x$ in $G$, the weight of $x$, $w(x) = \sum_{x\in V(G)} f(x) - (i+(v+1-i)) = \sum_{x\in V(G)} f(x) - (v+1)$, which is independent of the choice of $x$. Therefore, $G$ is $\{1,2,\ldots,d-1\}$-distance magic.
\end{proof}

\bl
Every connected graph is $\{0,1,\ldots,d\}$-distance magic.
\label{trivial}
\el
\begin{proof} The proof is straightforward since under the $\{0,1,\ldots,d\}$-distance magic labeling, we sum all labels in the graph in counting the weight of a vertex.
\end{proof}

\noi For obvious reason, we shall call the $\{0,1,\ldots,d\}$-distance magic of $G$ the {\em trivial $D$-distance magic labeling} of $G$. The following lemma deals with similar result for non-connected graphs.
\bl
Let $G$ be a non-connected graph having connected components $G_1,G_2, \ldots, G_p$, each of diameter $d_1, d_2, \ldots, d_p$, respectively. Let $d_{max}=\max_i d_i$ and $|V(G_i)|=n$ for each $i$. $G$ is $\{0,1,\ldots,d_{max}\}$-distance magic if and only if $n$ is even or both $n$ and $p$ are odd.
\label{disconnected:0-k}
\el
\begin{proof}
Suppose $G$ is $\{0,1,\ldots,d_{max}\}$-distance magic. Since the weight of a vertex $x$ is the sum of all labels in the component containing $x$, then such a sum must equal to the magic constant {\sf k}. Now we count the sum of all labels by two different ways of counting:
\begin{eqnarray}
{\sf k} p & = & 1 + \ldots + np \nonumber \\
{\sf k} p & = & \frac{(np+1)(np)}{2} \nonumber \\
{\sf k} & = & \frac{(np+1)n}{2}. \nonumber
\end{eqnarray}
To guarantee that both sides are integers then $n$ has to be even or both $n$ and $p$ must be odd.

\noi To prove the sufficiency, let $x_{ij}, 1\leq j \leq n$, be the vertices in the component $G_i$. If $n$ is even, label the vertices in the following way
\[f(x_{ij})= \left\{ \begin{array}{ll}
             i+(j-1)p,     & i \ {\rm odd,}\\
             p-i+1+(j-1)p, & i \ {\rm even.}\\
                     \end{array}
             \right.\]
With this labeling, the sum of all labels in the component $G_i$ is $\frac{n}{2}(np+1)$, which is equal to $w(x)$, for $x$ a vertex in $G_i$.

\noi If $n$ is odd, consider $n=2k+1$, $p=2m+1$, and the labeling $f$ as defined bellow.
\[f(x_{ij})= \left\{ \begin{array}{ll}
 2i-1,         & 1\leq i \leq m+1 \ {\rm and} \ j=1,\\
 2(i-m-1),     & m+2\leq i \leq 2m+1 \ {\rm and} \ j=1,\\
 4m+3-i,       & 1\leq i \leq 2m+1 \ {\rm and} \ j=2,\\
 5m+4-i        & 1\leq i \leq m+1  \ {\rm and} \ j=3,\\
 7m+5-i        & m+2\leq i \leq 2m+1  \ {\rm and} \ j=3,\\
 i+(j-1)(2m+1),& 1\leq i \leq 2m+1 \ {\rm and} \ j>3, \ j \ {\rm even,}\\
 2m+2-i+(j-1)(2m+1), & 1\leq i \leq 2m+1 \ {\rm and} \ j>3, \ j \ {\rm odd.}
                     \end{array}
             \right.\]
Thus, the sum of all labels in the component $G_i$ is $(9m+6)+(k-1)(2m+2)+(k-1)(2k+3)(2m+1)$.
\end{proof}

\noi In the previous lemma, we only consider graphs having all connected components of the same order. As to graphs having connected components with different order, we have $K_2 \cup K_1$ as an example of $\{0,1\}$-distance magic graph. Whether there are other graphs remains a question.
\bop
Let $G$ be a non-connected graph having connected components $G_1,G_2, \ldots, G_p$, each of diameter $d_1, d_2, \ldots, d_p$, respectively. Let $d_{max}=\max_i d_i$ and there exist $i,j$ such that $|V(G_i)|\neq |V(G_j)|$. Does there exist $G$ admitting $\{0,1,\ldots,d_{max}\}$-distance magic labeling other than $K_2 \cup K_1$?
\eop

\noi The last result in this section provides connection between $D$-distance magic labelings with different $D$s.
\bl {\em \cite{OS13}}
Let $D\subseteq \{0,1, \ldots, d\}$ and $D^\ast = \{0,1, \ldots, d\}-D$. Then $G$ is $D$-distance magic if and only if $G$ is $D^\ast$-distance magic.
\label{complement}
\el

\noi As a consequence of Lemma \ref{complement}, we have the following.
\bl
A graph of diameter $d$ is not $\{1,2,\ldots ,d\}$-distance magic.
\label{complement0}
\el

\noi We shall call the $D^\ast$-distance magic labeling in Lemma \ref{complement} the {\em complement labeling} of $D$-distance magic labeling. In the  following we extend the result to non-connected graphs.
\bl
Let $G$ be a graph having connected components $G_1,G_2, \ldots G_p$ of diameters $d_1, d_2, \ldots d_p$, respectively. Let $D\subseteq \{0,1, \ldots, d_{max}\}$ and $D^\ast = \{0,1, \ldots, d_{max}\}-D$, where $d_{max}=\max_i d_i$. If $G$ admits a $D$-distance magic labeling $f$ such that  $\sum_{x\in G_i} f(x)$ is constant for each $i$ then $G$ is $D^\ast$-distance magic. Conversely, if $G$ admits a $D^\ast$-distance magic labeling $f^\ast$ such that  $\sum_{x\in G_i} f^\ast(x)$ is constant for each $i$ then $G$ is $D$-distance magic.
\label{disconnectedcomplement}
\el
\begin{proof} For each $x\in V(G)$, we define $w(x)=\sum_{u\in N_D(x)} f(u)$ and $w^\ast(x)=\sum_{u\in N_{D^\ast}(x)} f(u)$. Clearly $w^\ast(x)=\sum_{u\in G_{x}} f(u) - w(x)$, where $G_{x}$ is the component containing $x$. If $w(x)$ is constant for each vertex $x$, then so is $w^\ast(x)$. The converse can be proved similarly.
\end{proof}

\noi In the next section, we study the existence of $D$-distance magic labelings with various $D$ for cycles.

\section{$D$-distance magic labelings for cycles}

\noi We shall start with cycles of even order.
\bt
Every even cycle $C_{2k}$ is $\{1,2,\ldots,k-1\}$-distance magic.
\label{C2k}
\et
\begin{proof} 
Each vertex in $C_{2k}$ is at distance $k$ from exactly one other vertex and so $C_{2k}$ is $\{1,2,\ldots,k-1\}$-distance magic by Lemma \ref{1k-1}.
\end{proof}
\noi As a direct consequence of Lemma \ref{complement}, we obtain
\bc
Every even cycle $C_{2k}$ is $\{0,k\}$-distance magic.
\label{C2kcomp}
\ec

\noi The next result is a characterization of cycles admitting $D$-distance magic labelings where $D$ is a singleton. 
\bt
For $k$ a positive integer, a cycle $C_n$ is $\{k\}$-distance magic if and only if $n=4k$.
\label{C4k}
\et
\begin{proof} 
Suppose that $f$ is a $\{k\}$-distance magic labeling of $C_n$. Let $x$ be an arbitrary vertex in $C_n$, then there exist exactly two vertices of distance $k$ from $x$, say $x_1$ and $x_2$. There also exists another vertex of distance $k$ from $x_1$ beside $x$, say $y$, and similarly there exists another vertex of distance $k$ from $y$ beside $x_1$, say $y_2$. If $n\neq 4k$ then $x, x_1, x_2, y, y_2$ are all distinct. Thus we obtain a contradiction by Lemma \ref{d-1}.

\noi To proof the sufficiency, suppose that $n=4k$. Notice that each vertex $x$ in $C_{4k}$ has a distinct pair of vertices of distance $k$, say $x_1$ and $x_2$. We label such a pair with a labeling $f$ such that $f(x_1) = 4k + 1 - f(x_2)$. Thus the weight of $x$, $w(x)=f(x_1)+f(x_2)=4k$ (independent of the choice of $x$) and so $C_{4k}$ is $\{k\}$-distance magic.
\end{proof}
\noi As a direct consequence of Lemma \ref{complement}, we obtain
\bc
For $k$ a positive integer, a cycle $C_n$ is $\{0,1, \ldots, k-1, k+1, \ldots, \lfloor\frac{n}{2}\rfloor\}$-distance magic if and only if $n=4k$.
\label{C4kcomp}
\ec

\noi We could generalize the result in Theorem \ref{C4k} to 2-regular graphs which is a generalization of a result in \cite{OS11a}.
\bt \emph{\cite{OS11a}}
A 2-regular graph is $\{1\}$-distance magic if and only if it is the union of 4-cycles.
\et
\bt
For $k$ a positive integer, a 2-regular graph is $\{k\}$-distance magic if and only if it is a disjoint union of $C_{4k}$s.
\label{mC4k}
\et
\begin{proof}
The proof is similar to that of Theorem \ref{C4k}, except for proving the sufficiency, where we use the labeling $f$ such that $f(x_1) = m4k + 1 - f(x_2)$, where $m$ is the number of copies of $C_{4k}$.
\end{proof}

%
\noi Some additional negative results for cycles are presented in the following theorem and corollary. The next result is proved by using Lemma \ref{d-1}.
\bt
For $n \geq 2k+2$, a cycle $C_{n}$ is not $\{0,1, \ldots, k\}$-distance magic.
\et
\noi By Lemma \ref{complement}, we obtain
\bc
For $n \geq 2k+2$, a cycle $C_{n}$ is not $\{k+1, k+2, \ldots, \lfloor\frac{n}{2}\rfloor\}$-distance magic.
\ec

\noi We then have the problem of characterizing $D$-distance magic cycles, or more generally, $D$-distance magic 2-regular graphs.
\bop
Given a particular distance set $D$, what are the necessary and sufficient conditions for 2-regular graphs to have $D$-distance magic labeling?
\label{2reg}
\eop

\section{$D$-distance magic labelings with $D\subseteq \{0,1,2\}$}

\noi A well-known result of Blass and Harary \cite{BH79} stated that almost all graphs have diameter 2. Therefore in this section we dedicate our study to $D$-distance magic labelings where $D\subseteq \{0,1,2\}$. Since $\{1\}$-distance magic and $\{0,1\}$-distance magic labelings have been studied extensively, we only consider $D\in \{\{2\}, \{0,2\}, \{1,2\}$, $\{0,1,2\}\}$.

\noi In the next lemma, we shall present necessary conditions of the existence of $D$-distance magic graphs with $D$ containing 2 but not 0.
\bl
Let $D$ be a distance set containing 2 but not 0. If $G$ is a graph of diameter at least 2 containing either
\begin{enumerate}
  \item two adjacent pendants, or
  \item two vertices of distance 2 having the same neighborhood,
\end{enumerate}
then $G$ is not $D$-distance magic.
\label{not-D}
\el
\begin{proof} 
Suppose $G$ is $D$-distance magic and let $x,y$ be the two adjacent pendants (in case 1) or the two vertices of distance 2 having the same neighborhood (in case 2). In both cases, since $D$ containing 2 but not 0, $N_D(x)$ and $N_D(y)$ containing exactly the same vertices except $x$ for $N_D(x)$ and $y$ for $N_D(y)$. Thus since $w(x)=w(y)$, we have $f(x)=f(y)$, a contradiction.
\end{proof}

\noi By the aforementioned lemma, many trees do not have $D$-labelings, where $D$ containing 2 but not 0. However, to characterize trees admitting such a labeling needs further study. More specifically, it is interesting to determine which trees have $D$-distance magic labelings where $D\subseteq \{0,1,2\}$.
\bop
What are the necessary and sufficient conditions for trees to have  $D$-distance magic labelings where $D\subseteq \{0,1,2\}$?
\eop

\subsection{$\{2\}$-distance magic labelings}

\bt
A complete multipartite graph is not $\{2\}$-distance magic.
\label{2Kmn}
\et
\begin{proof}
Let $x$ and $y$ be two vertices in the same partite set of a multipartite graph $G$. If we name the partite set $V_0$ then $N_{\{2\}}(x)=V_0 - \{x\}$ and $N_{\{2\}}(y)=V_0 - \{y\}$. By Lemma \ref{d-1}, $G$ is not $\{2\}$-distance magic.
\end{proof}

\noi Based on this result and the results of O'Neal and Slater \cite{OS11b} on extremal graphs of diameter 2 and 3, we suspect that graphs with diameter $2$ are not $\{2\}$-distance magic and more generally, graphs with diameter $d$ are not $\{d\}$-distance magic.
\bcon
Graphs with diameter $2$ are not $\{2\}$-distance magic. More generally, graphs with diameter $d$ are not $\{d\}$-distance magic.
\econ

\subsection{$\{0,2\}$-distance magic labelings}

\noi By Lemma \ref{complement}, we have the following results as consequences of the existence of $\{1\}$-distance magic labelings for particular graphs of diameter 2.

\bt {\em \cite{MRS03}} Let $H_{n,p}, n> 1 \ {\rm and} \ p> 1$, denote the complete symmetric multipartite graph with $p$ parts, each of which
contains $n$ vertices. $H_{n,p}$ is $\{0,2\}$-distance magic if and only if either $n$ is even or both $n$ and $p$ are odd.
\et

\bt {\em \cite{MRS03}}
Let $1\leq a_1 \leq a_2 \leq a_3$. Let $s_i = \sum_{j=1}^{i} a_j$, $p=2$ for complete bipartite graph $K_{a_1,a_2}$, and $p=3$ for complete tripartite graph $K_{a_1,a_2,a_3}$. There exist $\{0,2\}$-distance magic labelings for  $K_{a_1,a_2}$ and  $K_{a_1,a_2,a_3}$ if and only if the following conditions hold.
\begin{description}
\item[(a)] $a_2 \geq 2$,
\item[(b)] $v(v+1) \equiv 0 \bmod 2p$, and
\item[(c)] $\sum_{j=1}^{s_i} (n+1-j) \geq \frac{iv(v+1)}{2p}$ for $1\leq i
\leq p$.
\end{description}
\et

\bt {\em \cite{FKK06}}
An odd order $r$-regular graph of diameter 2 is $\{0,2\}$-distance magic if and only if $r$ is even and $2 \leq r \leq n-2$.
\et

\bt {\em \cite{FKK06}}
Let $G$ be an odd order regular graph of diameter 2 and $n$ be an odd positive integer. Then the graph
$G[\overline{K_n}]$ is $\{0,2\}$-distance magic.
\et

\noi The aforementioned theorem deal with odd order $G$; however for even order $G$, we have an example in which the composition of $G$ with $\overline{K_n}$ does have a $\{0,2\}$-distance magic labeling.
\bt {\em \cite{SAS09}}
For $n\geq 1, C_4[\overline{K_n}]$ is $\{0,2\}$-distance magic.
\et

\noi For graphs of diameter other than 2, we have the path of order 4, which is of diameter 3, admitting a $\{0,2\}$-distance magic labeling. This leads to the following question.
\bop
Does there exist a graph of diameter larger than 2, other than $P_4$, admitting $\{0,2\}$-distance magic labeling?
\eop

\subsection{$\{1,2\}$-distance magic labelings}

\noi By Lemma \ref{complement0}, $\{1,2\}$-distance magic labelings do not exist for graphs of diameter 2, and so in the following theorem we construct $\{1,2\}$-distance magic labelings for infinite families of graphs with diameter larger than 2.

\bt
There exists an infinite family of regular graphs with diameter 3 admitting $\{1,2\}$-distance magic labeling.
\et
\begin{proof} We construct a graph $G$ with $V(G)=\{x, x_1, x_2, \ldots , x_n, y, y_1, y_2, \ldots , y_n\}$ and $E(G)=\{xx_i, yy_i|1\leq i\leq n\}\cup \{x_iy_j|1\leq i\leq n, 1\leq j\leq n, i\neq j\}$. We can see that $G$ is an $n$-regular graph of order $2n+2$ and diameter 3. Moreover, each vertex has a unique vertex of distance 3: $y$ for $x$ and $y_i$ for $x_i$, $1\leq i\leq n$. By Lemma \ref{1k-1}, $G$ is $\{1,2\}$-distance magic.
\end{proof}

\noi The existence of non-regular $\{1,2\}$-distance magic graphs or $\{1,2\}$-distance magic graphs of larger diameter remain open as stated in the following.
\bop
Does there exists an infinite family of non-regular graphs admitting $\{1,2\}$-distance magic labeling?
\eop
\bop
Does there exists an infinite family of graphs with diameter at least 4 admitting $\{1,2\}$-distance magic labeling?
\eop

\subsection{$\{0,1,2\}$-distance magic labelings}

\noi By Theorem \ref{trivial}, every graph of diameter 2 admits the trivial $D$-distance magic labeling, i.e., an $\{0,1,2\}$-distance magic labeling. We could not find $\{0,1,2\}$-distance magic graphs of larger diameter, and so we ask the following question.

\bop
Does there exist a graph of diameter at least 3 admitting an $\{0,1,2\}$-distance magic labeling? 
\eop

\bibliographystyle{alpha}

\begin{thebibliography}{99}

\bibitem{ARSP04} B.D. Acharya, S.B. Rao, T. Singh and V. Parameswaran, Neighborhood magic graphs, {\em Proc. Nat. Conf. Graph Theory Combin. Algorithm} (2004).
    
\bibitem{AC} Marcin Anholcer and Sylwia Cichacz, Note on distance magic products $G \circ C_4$, preprint.

\bibitem{ACPT} Marcin Anholcer, Sylwia Cichacz, Iztok Peterin, and Aleksandra Tepeh, Distance magic labeling and two products of graphs, preprint.

\bibitem{AFK13} S. Arumugam, Dalibor Froncek, and N. Kamatchi, Distance Magic Graphs - A Survey, \emph{J. Indones. Math. Soc.} \textbf{Special Edition} (2011) 11-26.

\bibitem{Be09} S. Beena, On $\Sigma$ and $\Sigma ^{'}$ labelled graphs, {\em Disc. Math.}, {\bf 309} (2009) 1783-1787.

\bibitem{BH79} A. Blass, F. Harary, Properties of almost all graphs and complexes, {\em J. Graph Theory} {\bf 3} (1979) 225–240.

\bibitem{CF} Sylwia Cichacz and Dalibor Froncek, Distance magic circulant graphs, preprint.

\bibitem{CG} Sylwia Cichacz and Agnieszka Gorlich, Constant sum partition of sets of integers and distance magic graphs, preprint.

\bibitem{Fr07} D. Froncek, Fair incomplete tournaments with odd number of teams and large number of games, {\em Congressus Numerantium} \textbf{187} (2007) 83–89.

\bibitem{FKK06} D. Froncek, P. Kovar and T. Kovarova, Fair incomplete tournaments, {\em Bull. Inst. Combin. App.} {\bf 48} (2006) 31-33.

\bibitem{FKK11} D. Froncek, P. Kovar and T. Kovarova, Constructing distance magic graphs from regular graphs, \emph{J. Combin. Math. Combin. Comput.} \textbf{78} (2011) 349-354.

\bibitem{Ga} J. Gallian, A Dynamic Survey of Graph Labeling, {\em Electronic J. Combin.} \textbf{19} (2012) \#DS6.

\bibitem{GY92} J.R. Griggs and R.K. Yeh, Labelling graphs with a condition at distance 2, {\em SIAM J. Disc. Math.} {\bf 4} (1992) 586-595.

\bibitem{Ji99} M.I. Jinnah, On $\Sigma$-labelled graphs, {\em Technical Proceedings of Group Discussion on Graph Labeling Problems}, (1999)
71-77.

\bibitem{KR70} A. Kotzig and A. Rosa, Magic valuations of finite graphs, {\em Canad. Math. Bull.} {\bf 13} (1970) 451-461.

\bibitem{KFK12} P. Kovar, D. Froncek, and T. Kovarova,  A note on 4-regular distance magic graphs, \emph{Australasian J. Combin. }\textbf{54} (2012) 127-132.
    
\bibitem{KS12} Petr Kovar and Adam Silber, Distance magic graphs of high regularity, \emph{AKCE Int. J. Graphs Combin.} \textbf{9} (2012) 213-219.

\bibitem{MRS03} M. Miller, C. Rodger, and R. Simanjuntak, Distance magic labelings of graphs, {\em Australasian J. Combin.} {\bf 28} (2003) 305 - 315.

\bibitem{OS11a} A. O'Neal and P. Slater, The Minimax, Maximin, and Spread Values For Open Neighborhood Sums for 2-Regular Graphs, \emph{Math. Comput. Sci.} \textbf{5} (2011) 69–80.

\bibitem{OS11b} A. O'Neal and P. Slater, An introduction to distance $D$ magic graphs, \emph{J. Indonesian Math. Soc.} \textbf{Special Edition} (2011) 89-107.

\bibitem{OS13} A. O'Neal and P. Slater, Uniqueness Of Vertex Magic Constants, \emph{SIAM J. Disc. Math.} \textbf{27} (2013) 708–716

\bibitem{RSP04} S.B. Rao, T. Singh and V. Parameswaran, Some sigma labelled graphs I, {\em Graphs, Combinatorics, Algorithms and Applications} (2004) 125-133.

\bibitem{Ra08} S.B. Rao, Sigma Graphs - A survey, {\em Labelings of Discrete Structures and Applications} (2008) 135-140.

\bibitem{Sed64} J. Sedl\'{a}\v{c}ek, {\bf Problem 27} in Theory of Graphs and its Applications, {\em Proc. Symposium Smolenice 1963}, Prague (1964)
163-164.

\bibitem{SAA09} M. Seoud, A. E. I. Abdel Maqsoud and Y. I. Aldiban, New classes of graphs with and
without 1-vertex magic vertex labeling, {\em Proc. Pakistan Acad. Sci.} {\bf 46} (2009) 159-174.

\bibitem{SAS09} M.K. Shafiq, G. Ali and R. Simanjuntak, Distance magic labelings of a union
of graphs, {\em AKCE J. Graphs. Combin.} {\bf 6} (2009) 191-200.

\bibitem{SFMRW09} K.A. Sugeng, D. Froncek, M. Miller, J. Ryan and J. Walker, On distance magic
labeling of graphs, {\em J. Combin. Math. Combin. Comput.} {\bf 71} (2009) 39-48.

\bibitem{Vi94} V. Vilfred, {\em Sigma labelled graphs and circulant graphs}, Ph.D. Thesis, University of Kerala,
India (1994).

\end{thebibliography}

\end{document}